\documentclass[12pt]{article}
\usepackage{amsfonts}
\usepackage{latexsym, amssymb, amsmath, amscd, amsfonts, epsfig, graphicx, colordvi}
\usepackage{ifpdf}
\usepackage{graphicx}
\usepackage{pstricks}
\newtheorem{thm}{Theorem}[section]

\newtheorem{core}[thm]{Corollary}

\def\qed{\nopagebreak\hfill{\rule{4pt}{7pt}}
\medbreak}

\def\qed{\nopagebreak\hfill{\rule{4pt}{7pt}}
\medbreak}

\begin{document}

\begin{center}
{\bf \large Euler's Partition Theorem with
Upper Bounds

on Multiplicities }

\vskip 3mm
\end{center}

\begin{center}
William Y. C. Chen$^1$, Ae Ja Yee$^2$, and  Albert J. W. Zhu$^3$

 $^{1,3}$Center for Combinatorics, LPMC-TJKLC\\
   Nankai University\\
    Tianjin 300071, P.R. China

\vskip 1mm

$^2$Department of Mathematics\\
The Pennsylvania State University\\
University Park, PA 61802 USA

\vskip 1mm

   Email: $^1$chen@nankai.edu.cn, $^2$yee@math.psu.edu, $^3$zjw@cfc.nankai.edu.cn

\end{center}

\vskip 6mm \noindent{{\bf Abstract.}
We show that the number of partitions of $n$
with  alternating sum $k$ such that the multiplicity
of each part is
bounded by $2m+1$ equals the number of
 partitions of $n$  with $k$ odd parts
 such that the multiplicity
of each even part is bounded by $m$.
The first proof relies on
two formulas with two parameters that are
related to the four-parameter formulas of Boulet.
 We also give
a combinatorial proof of this result by using Sylvester's
bijection, which implies a stronger partition theorem.
 For $m=0$, our result reduces to
 Bessenrodt's refinement of Euler's Theorem.
If the alternating sum and the number of odd parts
are not taken into account, we are led to a connection
to a generalization of Euler's theorem,
which can be deduced from a theorem of Andrews on
 equivalent upper bound sequences of multiplicities.
Analogously, we show that the number of partitions of $n$ with
alternating sum $k$ such that the multiplicity
of each even part is bounded by $2m+1$ equals
the number of partitions of $n$ with
 $k$ odd parts such that the multiplicity
of each even part is also bounded by $2m+1$. {We provide a combinatorial proof as well.}



\section{Introduction}

The main objective of this paper is to present a unification
of a refinement and a generalization of Euler's partition
theorem. More precisely, we show that the number of partitions of
$n$ with a given alternating sum $k$
such that the multiplicity
of each part is bounded by $2m+1$
 is equal to the number of partitions of $n$
with a given number of odd parts $k$
 such that the multiplicity
of each even part is bounded by $m$.
 The algebraic proof of this result
relies on two formulas with two parameters
that are related to the four-parameter
formulas of Boulet
on the generating function for
partitions \cite{Bou06}. Meanwhile,
 we give a combinatorial proof by using
 Sylvester's bijection \cite{Syl82} which implies a
 stronger partition theorem.

For the case $m=0$, our result reduces
to a refinement of Euler's partition theorem
due to  Bessenrodt \cite{Bes94}.
From a different perspective,
our result is also related to a theorem of
Andrews \cite{And72} on equivalent
upper bound sequences of multiplicities. If the alternating
sum and the number of odd parts are not taken into account
in our theorem, we obtain a connection between
partitions with the multiplicity
of each part bounded by $2m+1$
and partitions with the multiplicity
of each even parts bounded by
$m$.  This is clearly a generalization of Euler's
theorem. Indeed, it is a special case of
a theorem of  Andrews. In the language of Andrews,
we may say that the upper bound sequence of multiplicities
 $(2m+2,\, 2m+2,\, 2m+2,\, \ldots)$ is equivalent to
 the upper bound sequence of  multiplicities
 $(\infty,\, m+1,\, \infty,\, m+1,\, \ldots)$.
 However, when the alternating
 sums and the number of odd parts
 are taken into consideration, there does not seem to
 be a general theorem on the equivalent upper
 bound sequences of multiplicities.

Using the same idea to prove our main theorem,
we  show that the number of
partitions of $n$ with alternating sum $k$
such that the multiplicity
of each even part is bounded by $2m+1$ equals
the number of partitions of $n$ with
 $k$ odd parts such that the multiplicity
of each even part is also bounded by $2m+1$.
{We provide a combinatorial proof as well.}
In particular, for $m=0$, we find that
the number of partitions of $n$ with alternating sum
$k$ such that the even parts are distinct
 equals the number of partitions of $n$ with $k$ odd parts
 such that the even parts are distinct.

Let us give an overview of the theorem
of Bessenrodt and the theorem of Andrews.
Throughout this paper,
 we shall adopt the common
  notation on partitions used in
Andrews \cite{And76}.
A partition $\lambda$ of a positive integer
$n$ is a
finite sequence of positive integers
$\lambda_1, \lambda_2, \ldots, \lambda_r$ such that
$\lambda_1\ge \lambda_2\ge \cdots \ge \lambda_r$ and
$\lambda_1+\lambda_2+\cdots +\lambda_r=n$.
We write $\lambda=(\lambda_1, \lambda_2, \ldots, \lambda_r)$.
Each $\lambda_i$ is called a part.
The part $\lambda_1$
 is called the largest part of $\lambda.$ The
number of parts of $\lambda$ is called the length
of $\lambda,$
denoted $l(\lambda).$
The weight of $\lambda$ is the sum of
parts of $\lambda$,  denoted $|\lambda|$.
 The conjugate partition
of $\lambda$ is defined by
$\lambda'=(\lambda'_1,\,\lambda'_2,\,\ldots,\,\lambda'_t),$
 where
$\lambda'_i$
is the number of parts of $\lambda$
greater than or equal to $i.$
The  number of
odd parts in $\lambda$ is denoted by  $l_o(\lambda)$.
The alternating sum of $\lambda,$ namely,
$$l_a(\lambda)=
\lambda_1-\lambda_2+\lambda_3-\lambda_4+\cdots.$$
is denoted by $l_a(\lambda)$.

Euler's partition theorem states that
the number of partitions of $n$ into
distinct parts is equal to the number of
partitions of $n$ into odd parts.
There are many refinements
of Euler's partition theorem, see, for examples,
Bessenrodt \cite{Bes94}, Fine \cite[pp.46-47]{Fin88}
and Sylvester
\cite[p.24]{And76}. For $m=0$,
 our partition
theorem reduces
to a refinement of Euler's theorem due to Bessenrodt
\cite{Bes94}, which is a limiting
 case of the lecture hall theorem due
  to  Bousquet-M\'{e}lou
and Erikssonin \cite{Bou971,Bou972}.
Bessenrodt derived the following relation  from
Sylvester's bijection \cite{Syl82}.

\begin{thm} \label{euler-bousquet}
The number of partitions of $n$ into distinct parts
with alternating sum $l$ is equal to the number of
partitions of $n$ into odd parts with length $l$.
 In terms of generating functions, we have
\begin{equation}\label{euler-bousquet-e}
\sum_{\lambda \in \mathcal{D}}
y^{l_a(\lambda)}q^{|\lambda|}=\sum_{\mu \in
\mathcal{O}}y^{l(\mu)}q^{|\mu|}.
\end{equation}
\end{thm}

There is another refinement of Euler's theorem due to
Glaisher \cite{Gla83}.
A unification of the refinements of
Bessenrodt \cite{Bes94} and Glaisher \cite{Gla83}
 has been obtained by Chen, Gao, Ji
and Li \cite{Chen10} which can be deduced from
Boulet's four-parameter formula \cite{Bou06}.
A combinatorial proof of the unification
is given in \cite{Chen10}.

In another direction, several
generalizations of Euler's theorem have been obtained,
see,
for example,
Alder \cite{Ald69}, Andrews \cite{And66,And72,AndE04},
Glaisher \cite{Gla83} and Moore \cite{Moo74}.
Andrews \cite{And72} proved the following theorem on
equivalent upper bound sequences of multiplicities, see also Pak
\cite{Pak06}.
An upper bound sequence of multiplicities is
defined to be an infinite sequence
 $a=(a_1,a_2,a_3,\ldots)$ of nonnegative
integers or infinity.
 Let $supp(a)$ be the set of indices $i$
 such that $a_i$ is finite.  We say that two upper bound sequences of multiplicities
  $a=(a_1,a_2,a_3,\ldots)$ and $b=(b_1,b_2,b_3,\ldots)$
  are equivalent,
  if there exists a one-to-one correspondence
 $\pi\colon supp(a)\rightarrow supp(b)$
 such that $i a_i=j b_j,$
 for all $j=\pi(i).$

\begin{thm}
 \label{andrews-partition}
 Let
  $a=(a_1,a_2,a_3,\ldots)$ and
 $b=(b_1,b_2,b_3,\ldots)$ be two
 upper bound sequences of multiplicities, and
  let $\mathcal{A}_n$ and $\mathcal{B}_n$ be the sets
 of partitions $\lambda=(1^{m_1}2^{m_2}\cdots)$ of $n$
 such that $m_i<a_i$
 and $m_i<b_i$ for all $i$. If $a$ and $b$ are equivalent,
  then
 $|\mathcal{A}_n|=|\mathcal{B}_n|,$ for all $n>0.$
\end{thm}

In the language of Andrews, we may say that
the upper bound sequence of multiplicities
 $a=(2m+2,\, 2m+2,\, 2m+2,\, \ldots)$ is equivalent
 to the upper bound sequence of multiplicities
 $b=(\infty,\, m+1,\, \infty,\, m+1,\, \ldots)$.
  To be more specific,
 let $\pi$ be the map $i\rightarrow 2i$,
 then we have  $i a_i=\pi(i) b_{\pi(i)}$.  In this case,
 we have the following consequence, which can be
 considered as a generalization of Euler's theorem.

\begin{core}\label{core of andrews}
The number of
 partitions of $n$ with the multiplicity of
 each part bounded by $2m+1$
equals the number of partitions of $n$
with the multiplicity of each even part bounded by $m.$
\end{core}

It can be seen that our main result is   a unification
of Theorem \ref{euler-bousquet} and
Corollary \ref{core of andrews}.
This paper is organized as follows. In Section 2,
we provide two formulas on two-parameter generating
functions for partitions that are related to Boulet's
four-parameter formula by imposing an upper bound
sequence of multiplicities.
In Section 3, we use the two formulas presented
in Section 2 to prove the main theorem.
We also prove the second result of this paper,
that is, with the restrictions of an upper
bound $2m+1$ on the multiplicity of
each even part,
the number of partitions of $n$ with
alternating sum $k$ equals the number of partitions of $n$
with $k$ odd parts.
In Section 4, we give a combinatorial proof of the
main result based on Sylvester's bijection and
obtain a stronger version concerning the number of odd parts,
the largest odd part and the largest part with odd multiplicity. {Finally, we provide a combinatorial proof of our second theorem in the last section.}

\section{Two-parameter formulas}

To prove the main result of this paper, we need
two formulas on  two-parameter generating functions
for partitions. The first formula (Theorem \ref{3-1})
can be deduced from the Boulet's
 four-parameter formula \cite{Bou06}.
It seems that there are
typos in the formula of Boulet.
A corrected version with a technical condition added
is presented in Theorem \ref{boulet-general}.
The first two-parameter formula (Theorem \ref{3-1})
can be deduced from Theorem \ref{boulet-general}
for the case $i=0$, $k=1$, $a=b$ and $c=d$.
We provide the second  two-parameter formula
(Theorem \ref{general two parameter})
from a two-parameter formula which is the case $a=c$ and $b=d$
of (\ref{Boulet-four-parameter}). The proof is similar to that of
Boulet's four-parameter formula \cite{Bou06}.
It should be noted that Theorem \ref{3-1} can be
considered as an extension of the two-parameter
formula obtained from (\ref{Boulet-four-parameter}) by setting $a=b$ and $c=d$.
The proof of the first two-parameter formula (Theorem \ref{3-1}) is analogous to that of
Theorem \ref{general two parameter}, and hence it is omitted.

To state Boulet's four-parameter formula, let $\lambda$ be a partition,
and let $a,$ $b,$ $c,$ $d$ be commuting indeterminants.
Define
\begin{equation}\label{four-parameter}
\omega(\lambda)=a^{\sum_{i\ge 1}\lceil{\lambda_{2i-1}/2}\rceil}
b^{\sum_{i\ge 1}\lfloor{\lambda_{2i-1}/2}\rfloor}
c^{\sum_{i\ge 1}\lceil{\lambda_{2i}/2}\rceil}
d^{\sum_{i\ge 1}\lfloor{\lambda_{2i}/2}\rfloor}.
\end{equation}
The
above four-parameter weight is introduced
by Boulet \cite{Bou06}.
Let \[ \Phi(a,b,c,d)=\sum_\lambda \omega(\lambda),\]
where the sum
is over all partitions. Boulet \cite{Bou06}
obtained the
following  formula.

\begin{thm} \label{Boulet} We have
\begin{equation}\label{Boulet-four-parameter}
\Phi(a,b,c,d)=\prod_{j=1}^{\infty}
\frac{(1+a^{j}b^{j-1}c^{j-1}d^{j-1})
(1+a^{j}b^{j}c^{j}d^{j-1})}
{(1-a^{j}b^{j}c^{j}d^{j})
(1-a^{j}b^{j}c^{j-1}d^{j-1})
(1-a^{j}b^{j-1}c^{j}d^{j-1})}.
\end{equation}
\end{thm}

Boulet generalized the above formula by
consider an upper bound sequence of
multiplicities. We shall
adopt the notion in \cite{Bou06}.
Let $R$ be a subset of positive integers,
and let $\rho$ be a map from $R$ to the set $E$ of even
 positive integers. Boulet defined
 Par$(i,k;R,\rho)$ as the set
 of  partitions $\lambda$ with
parts congruent to $i$ $\mbox{(mod $k$)}$ such
that for any $r\in R$  the multiplicity of $r$ in $\lambda$
is less than $\rho(r)$. In other words, $\rho$ can be
viewed as an upper bound sequence of multiplicities for
partitions in Par$(i,k;R,\rho)$. Let
\begin{equation}\label{2-1}
\Phi_{i,k;R,\rho}(a,b,c,d)
=\sum_{\lambda}\omega(\lambda),
\end{equation}
where the sum is over all partitions
in Par$(i,k;R,\rho)$.

Boulet \cite{Bou06} derived a formula for
$\Phi_{i,k;R,\rho}(a,b,c,d)$, which seems to contain
typos. This formula  can be corrected
as follows. First, we impose a further  condition
on the definition of Par$(i,k;R,\rho)$, that is,
for  $i\not=0$, the length of $\lambda$ is
even and the part $i$ appears at most once. Then the
formula for $\Phi_{i,k;R,\rho}(a,b,c,d)$ takes
the following form.

\begin{thm}\label{boulet-general} We have
\begin{equation}\label{2-2}
\Phi_{i,k;R,\rho}(a,b,c,d)=ST,
\end{equation}
where
$$S=\prod_{j=1}^{\infty}
\frac{(1+a^{\lceil\frac{jk+i}{2}\rceil}
b^{\lfloor\frac{jk+i}{2}\rfloor}
c^{\lceil\frac{(j-1)k+i}{2}\rceil}
d^{\lfloor\frac{(j-1)k+i}{2}\rfloor})}
{(1-a^{\lceil\frac{jk+i}{2}\rceil}
b^{\lfloor\frac{jk+i}{2}\rfloor}
c^{\lceil\frac{jk+i}{2}\rceil}
d^{\lfloor\frac{jk+i}{2}\rfloor})
(1-a^{jk}b^{jk}c^{(j-1)k}d^{(j-1)k})}$$
and
$$T=\prod_{r\in R}(1-a^{\lceil\frac{r}{2}\rceil
\frac{\rho(r)}{2}}
b^{\lfloor\frac{r}{2}\rfloor\frac{\rho(r)}{2}}
c^{\lceil\frac{r}{2}\rceil\frac{\rho(r)}{2}}
d^{\lfloor\frac{r}{2}\rfloor\frac{\rho(r)}{2}}).$$
\end{thm}
For the purpose of this paper, we only need
the special case of (\ref{2-2}), which is the case for $k=1$, $i=0$, $a=b$ and $c=d$.
For this special case,  we define
\begin{equation}\label{two-parameter-1}
\omega_1(\lambda)=a^{\sum_{i\ge 1}\lceil{\lambda_{2i-1}/2}\rceil}
a^{\sum_{i\ge 1}\lfloor{\lambda_{2i-1}/2}\rfloor}
b^{\sum_{i\ge 1}\lceil{\lambda_{2i}/2}\rceil}
b^{\sum_{i\ge 1}\lfloor{\lambda_{2i}/2}\rfloor}.
\end{equation}
For notational simplicity, we denote by Par$(R,\rho)$
 the set of  partitions such that any
 element $r$ in $R$ appears as
a part with multiplicity less than  $\rho(r)$.
Set
\[ \Phi_{R,\rho}(a,a,b,b)
=\Phi_{0,1;R,\rho}(a,a,b,b)=
\sum_{\lambda}\omega_1(\lambda),\]
 where the
sum is over all partitions in Par$(R,\rho).$
We can deduce the following formula from (\ref{2-2}).

\begin{thm}\label{3-1} We have
\begin{align*}
\Phi_{R,\rho}(a,a,b,b)
=&\prod_{j=1}^{\infty}\frac{(1+a^{\lceil\frac{j}{2}\rceil
+\lfloor\frac{j}{2}\rfloor}
b^{\lceil\frac{j-1}{2}\rceil
+\lfloor\frac{j-1}{2}\rfloor})}
{(1-a^{\lceil\frac{j}{2}\rceil
+\lfloor\frac{j}{2}\rfloor}
b^{\lceil\frac{j}{2}\rceil
+\lfloor\frac{j}{2}\rfloor})
(1-a^{2j}b^{2j-2})}\\
&\quad \times\prod_{r\in R}
(1-a^{(\lceil\frac{r}{2}\rceil+\lfloor\frac{r}{2}\rfloor)
\frac{\rho(r)}{2}}
b^{(\lceil\frac{r}{2}\rceil+\lfloor\frac{r}{2}\rfloor)
\frac{\rho(r)}{2}}).
\end{align*}
\end{thm}

Bear in mind that  the condition that $\rho(r)$ is even
for any $r\in R$ is required in Theorem \ref{3-1}.
However, as will  be  seen, we can obtain
similar formulas for any upper bound sequence
of multiplicities as long as we have $a=c$ and $b=d$.
So we define
\begin{equation}\label{two-parameter-2}
\omega_2(\lambda)=a^{\sum_{i\ge 1}\lceil{\lambda_{2i-1}/2}\rceil}
b^{\sum_{i\ge 1}\lfloor{\lambda_{2i-1}/2}\rfloor}
a^{\sum_{i\ge 1}\lceil{\lambda_{2i}/2}\rceil}
b^{\sum_{i\ge 1}\lfloor{\lambda_{2i}/2}\rfloor}.
\end{equation}
To be more specific, let $\gamma$ be a map
from $R$ to positive integers, and let
\[ \Psi_{R,\gamma}(a,b,a,b)=\sum_{\lambda}\omega_2(\lambda),\]
 where the sum is over   partitions in
 Par($R,\gamma$). To derive our main theorem,
  we   need the following extension of a two-parameter
formula, that is, the generating function for
$\Psi_{R,\gamma}(a,b,a,b).$  The proof of the following two-parameter
formula is analogous to that of Boulet's four-parameter formula  \cite{Bou06}.

\begin{thm}\label{general two parameter}
We have
$$\Psi_{R,\gamma}(a,b,a,b)=UV,$$
where
$$U=\prod_{j=1}^{\infty}
\frac{(1+a^{\lceil\frac{j}{2}\rceil}
b^{\lfloor\frac{j}{2}\rfloor}
a^{\lceil\frac{j-1}{2}\rceil}
b^{\lfloor\frac{j-1}{2}\rfloor})}
{(1-a^{2\lceil\frac{j}{2}\rceil}
b^{2\lfloor\frac{j}{2}\rfloor})
(1-a^{2j-1}b^{2j-1})}$$
and
$$V=\prod_{r\in R}(1-a^{\lceil\frac{r}{2}\rceil
{(\gamma(r)+1)}}
b^{\lfloor\frac{r}{2}\rfloor(\gamma(r)+1)}).$$
\end{thm}

{\noindent \it Proof.}
In view of the definition of  $\Psi_{R,\gamma}(a,b,c,d),$
 it is easy to see that
 $\Psi_{\emptyset,\gamma}(a,b,a,b)$ equals
$\Phi(a,b,a,b)$ as given in Theorem \ref{Boulet}.
 Hence we obtain
\begin{align*}
\Psi_{\emptyset,\gamma}(a,b,a,b)
&=\prod_{j=1}^{\infty}
\frac{(1+a^{j}b^{j-1}a^{j-1}b^{j-1})
(1+a^{j}b^{j}a^{j}b^{j-1})}
{(1-a^{j}b^{j}a^{j}b^{j})(1-a^{j}b^{j}a^{j-1}b^{j-1})
(1-a^{j}b^{j-1}a^{j}b^{j-1})}\\[5pt]
&=\prod_{j=1}^{\infty}
\frac{(1+a^{\lceil\frac{j}{2}\rceil}
b^{\lfloor\frac{j}{2}\rfloor}
a^{\lceil\frac{j-1}{2}\rceil}
b^{\lfloor\frac{j-1}{2}\rfloor})}
{(1-a^{2\lceil\frac{j}{2}\rceil}
b^{2\lfloor\frac{j}{2}\rfloor})
(1-a^{2j-1}b^{2j-1})},
\end{align*}
which is the generating function for Par($\emptyset,\gamma$).

In order to obtain the generating function for Par($R,\gamma$),
we notice that   any partition $\lambda$ in Par($\emptyset,\gamma$)
   can be expressed as a bipartition ($\mu$, $\nu$)
  such that $\mu\in \mbox{Par}(R,\gamma)$ and $\nu\in\mathcal {L}$
  where $\mathcal {L}$ is the
set of partitions with parts $r\in R$
occurring a multiple of $\gamma(r)+1$ times. In other words, there is a bijection between
Par($\emptyset,\gamma$) and
Par($R,\gamma$)$\times\mathcal {L}$.
Clearly, any part $r\in R$ occurring a multiple
of $\gamma(r)+1$ times can be represented as
a multiple of the following block
$$\left.\underbrace{\begin{array}{ccccccc}
a&b&a&b&\cdots &a&b\\
a&b&a&b&\cdots &a&b\\
\cdots&\cdots&\cdots&\cdots&\cdots&\cdots&\cdots\\
a&b&a&b&\cdots &a&b
\end{array}}_{\mbox{$r$}}\right\}
\quad\mbox{repeated $\gamma(r)+1$ times.}$$
The weight of the above block is
$a^{\lceil\frac{r}{2}\rceil{(\gamma(r)+1)}}
b^{\lfloor\frac{r}{2}\rfloor(\gamma(r)+1)}$.
Thus the generating function for $\mathcal{L}$ equals
$$V^{-1}=\prod_{r\in R}\frac{1}
{(1-a^{\lceil\frac{r}{2}\rceil{(\gamma(r)+1)}}
b^{\lfloor\frac{r}{2}\rfloor(\gamma(r)+1)})}$$
It follows that
$$\Psi_{R,\gamma}(a,b,c,d)=UV.$$
This completes the proof.  \qed

\section{Two Partition Theorems}

In this section, we provide two partition theorems. The first one is stated as follows.

\begin{thm}\label{partition theorem1}
For $k, m, n \ge 0$, the number of partitions of $n$ with alternating sum $k$ such that
 each part appears at most $2m+1$ times equals the
number of partitions of $n$ with $k$ odd parts such that each even part
 appears at most $m$ times.
\end{thm}

For example, the following table  illustrates the case of $n=7$
and $m=1$.
\begin{center}
    \begin{tabular}{|c|c||c|c|} \hline
     $l_a$ & Each part appears at most 3 times &  $l_o$ & Each even part appears only once \\\hline
    $1$& $(2^2,1^3)$ $(2^3,1)$ $(3,2,1^2)$ $(3^2,1)$ $(4,3)$     &$1$& $(4,2,1)$ $(4,3)$ $(5,2)$ $(6,1)$ $(7)$\\\hline
    $3$& $(3,2^2)$ $(4,1^3)$ $(4,2,1)$ $(5,2)$   &$3$& $(3,2,1^2)$ $(3^2,1)$ $(4,1^3)$ $(5,1^2)$\\\hline
    $5$& $(5,1^2)$ $(6,1)$    &$5$& $(2,1^5)$ $(3,1^4)$\\\hline
    $7$& $(7)$       &$7$& $(1^7)$\\\hline
    \end{tabular}
\end{center}

If the alternating sum and the number of odd parts are not
taken into account, we are led to
Corollary \ref{core of andrews}, which can be regarded as a
generalization of Euler's theorem and a special case
 of Theorem \ref{partition theorem1}.
Meanwhile, it is easy to see that Theorem \ref{partition theorem1}
reduces to Theorem \ref{euler-bousquet} when $m=0.$

The second partition theorem is concerned with the alternating sum,
the number of odd parts and  an upper
bound $2m+1$ on the multiplicity of each even part.

\begin{thm}\label{partition theorem2}
For $k, m, n\ge 0$,
the number of partitions of $n$ with alternating sum $k$
such that each even part appears at most $2m+1$ times equals
the number of partitions of $n$ with $k$ odd parts
such that each even part appears at most $2m+1$ times.
\end{thm}

For example, for $n=7$
and $m=0$, the partitions in Theorem \ref{partition theorem2} are listed in the following table.
\begin{center}
    \begin{tabular}{|c|c||c|c|} \hline
     $l_a$ & Each even part appears only once &  $l_o$ & Each even part appears only once \\\hline
    $1$& $(1^7)$ $(2,1^5)$ $(3,2,1^2)$ $(3^2,1)$ $(4,3)$ & $1$& $(4,2,1)$ $(4,3)$ $(5,2)$ $(6,1)$ $(7)$ \\\hline
    $3$& $(3,1^4)$ $(4,1^3)$ $(4,2,1)$ $(5,2)$  &  $3$& $(3,2,1^2)$ $(3^2,1)$ $(4,1^3)$ $(5,1^2)$\\\hline
    $5$& $(5,1^2)$ $(6,1)$  &  $5$& $(2,1^5)$ $(3,1^4)$\\\hline
    $7$& $(7)$  &  $7$& $(1^7)$\\\hline
    \end{tabular}
\end{center}

To give the proofs of the above theorems, we shall adopt  the common notation
for the $q$-shifted factorials
\begin{equation*}\label{q notation finite}
(a;q)_{n}=(1-a)(1-aq)\cdots(1-aq^{n-1}),\quad n\ge 1,
\end{equation*}
and
\begin{equation*}\label{q notation infinite}
(a;q)_{\infty}=\lim_{n\rightarrow\infty}
(1-a)(1-aq)(1-aq^2)\cdots,
\end{equation*}
where $0\le|q|<1,$
see Gasper and Rahman \cite{Gas90}.

Notice that the generating functions for   partitions
in Theorem \ref{partition theorem1} and Theorem \ref{partition theorem2}
can be expressed in the notation
$\Phi_{R,\rho}(a,b,c,d)$ and
$\Psi_{R,\gamma}(a,b,c,d)$ defined
in the previous section, where $\rho$ and $\gamma$ are upper bound sequences
of multiplicities.
More precisely,  $\Phi_{R,\rho}(a,b,c,d)$ is the generating function for partitions
such that each part $r\in R$ appears less than $\rho(r)$ times with $\rho(r)$ being an even number,
 and $\Psi_{R,\gamma}(a,b,c,d)$ is the generating function for partitions
such that each part $r\in R$ appears at most $\gamma(r)$ times.

In the above notation, it is easy to check that
$$\Phi_{\mathbb{N},2m+2}(xq,xq,x^{-1}q,x^{-1}q)$$
is the generating function for partitions
such that each part appears at most $2m+1$ times and the exponent of $x$
is  the alternating sum.
We also observe that
$$\Psi_{E,m}(xq,x^{-1}q,xq,x^{-1}q)$$
 is  the generating function for partitions
such that each even part appears at most $m$ times and
 the exponent of $x$ is the number of odd parts.

Based on the above interpretations of the functions
$\Phi_{\mathbb{N},2m+2}(xq,xq,x^{-1}q,x^{-1}q)$ and
$\Psi_{E,m}(xq,x^{-1}q,xq,x^{-1}q)$, one can deduce Theorem \ref{partition theorem1}
from the following relations.

\begin{thm}\label{eq-gf}
We have
\begin{equation*}
\Phi_{\mathbb{N},2m+2}(xq,xq,x^{-1}q,x^{-1}q)
=\Psi_{E,m}(xq,x^{-1}q,xq,x^{-1}q)
=\frac{(-xq;q^2)_{\infty}(q^{2m+2};q^{2m+2})_{\infty}}{(q^2;q^2)_{\infty}
(x^2q^2;q^4)_{\infty}}.\label{4eq1}
\end{equation*}
\end{thm}

\allowdisplaybreaks

 {\noindent\it Proof.}
Setting $a=xq$, $b=xq$, $c=x^{-1}q$,
 $d=x^{-1}q$, $R=\mathbb{N},$ and
$\rho=2m+2$ in Theorem \ref{3-1},   we obtain
\begin{align*}
&\Phi_{\mathbb{N},2m+2}(xq,xq,x^{-1}q,x^{-1}q)\\[5pt]
=&\prod_{j=1}^{\infty}
\frac{(1+(xq)^{\lceil\frac{j}{2}\rceil}
(xq)^{\lfloor\frac{j}{2}\rfloor}
(x^{-1}q)^{\lceil\frac{j-1}{2}\rceil}
(x^{-1}q)^{\lfloor\frac{j-1}{2}\rfloor})}
{(1-(xq)^{\lceil\frac{j}{2}\rceil}
(xq)^{\lfloor\frac{j}{2}\rfloor}
(x^{-1}q)^{\lceil\frac{j}{2}\rceil}
(x^{-1}q)^{\lfloor\frac{j}{2}\rfloor})
(1-(xq)^{j}(xq)^{j}(x^{-1}q)^{j-1}
(x^{-1}q)^{j-1})}\\[5pt]
&\times\prod_{r\in \mathbb{N}}
(1-(xq)^{\lceil\frac{r}{2}\rceil
\frac{2m+2}{2}} (xq)^{\lfloor\frac{r}{2}\rfloor
\frac{2m+2}{2}} (x^{-1}q)^{\lceil\frac{r}{2}\rceil
\frac{2m+2}{2}} (x^{-1}q)^{\lfloor\frac{r}{2}\rfloor
\frac{2m+2}{2}})\\[5pt]
=&\frac{\prod_{j=1}^{\infty}
(1+(xq)^{j}(xq)^{j}(x^{-1}q)^{j}(x^{-1}q)^{j-1})
\prod_{j=1}^{\infty}(1+(xq)^{j}(xq)^{j-1}
(x^{-1}q)^{j-1}(x^{-1}q)^{j-1})}
{\prod_{j=1}^{\infty}(1-(xq)^{\lceil
\frac{j}{2}\rceil}(xq)^{\lfloor\frac{j}{2}\rfloor}
(x^{-1}q)^{\lceil\frac{j}{2}\rceil}
(x^{-1}q)^{\lfloor\frac{j}{2}\rfloor})
\prod_{j=1}^{\infty}(1-x^2q^{4j-2})}\\[5pt]
&\times\prod_{r\in \mathbb{N}}
(1-(xq)^{\lceil\frac{r}{2}\rceil(m+1)}
(xq)^{\lfloor\frac{r}{2}\rfloor(m+1)}
(x^{-1}q)^{\lceil\frac{r}{2}\rceil(m+1)}
(x^{-1}q)^{\lfloor\frac{r}{2}\rfloor(m+1)})\\[5pt]
=&\frac{\prod_{j=1}^{\infty}(1+xq^{4j-1})
\prod_{j=1}^{\infty}(1+xq^{4j-3})}
{\prod_{j=1}^{\infty}(1-q^{4j})(1-q^{4j-2})
\prod_{j=1}^{\infty}(1-x^2q^{4j-2})}
\times\prod_{r=1}^{\infty}(1-q^{4r(m+1)})
(1-q^{(4r-2)(m+1)})\\[5pt]
=&\frac{(-xq;q^2)_{\infty}(q^{2m+2};q^{2m+2})_{\infty}}{(q^2;q^2)_{\infty}
(x^2q^2;q^4)_{\infty}}.
\end{align*}
Meanwhile, setting $a=xq$,
$b=x^{-1}q$, $R=E$ (the set of even positive integers)  and
$\gamma=m$ in Theorem \ref{general two parameter}, we find that
\begin{align*}
&\Psi_{E,m}(xq,x^{-1}q,xq,x^{-1}q)\\[5pt]
=&\prod_{j=1}^{\infty}
\frac{(1+(xq)^{\lceil\frac{j}{2}\rceil}
(x^{-1}q)^{\lfloor\frac{j}{2}\rfloor}
(xq)^{\lceil\frac{j-1}{2}\rceil}
(x^{-1}q)^{\lfloor\frac{j-1}{2}\rfloor})}
{(1-(xq)^{2\lceil\frac{j}{2}\rceil}
(x^{-1}q)^{2\lfloor\frac{j}{2}\rfloor})
(1-(xq)^{2j-1}(x^{-1}q)^{2j-1})}\\[5pt]
&\times\prod_{r\in E}
(1-(xq)^{\lceil\frac{r}{2}\rceil(m+1)}
(x^{-1}q)^{\lfloor\frac{r}{2}\rfloor(m+1)})\\[5pt]
=&\frac{\prod_{j=1}^{\infty}(1+(xq)^{j}
(x^{-1}q)^{j}(xq)^{j}(x^{-1}q)^{j-1})
\prod_{j=1}^{\infty}(1+(xq)^{j}
(x^{-1}q)^{j-1}(xq)^{j-1}(x^{-1}q)^{j-1})}
{\prod_{j=1}^{\infty}(1-(xq)^{2\lceil
\frac{j}{2}\rceil}(x^{-1}q)^{2\lfloor\frac{j}{2}\rfloor})
\prod_{j=1}^{\infty}(1-q^{4j-2})}\\[5pt]
&\times\prod_{r=1}^{\infty}(1-(xq)^{r(m+1)}
(x^{-1}q)^{r(m+1)})\\[5pt]
=&\frac{\prod_{j=1}^{\infty}(1+xq^{4j-1})
\prod_{j=1}^{\infty}(1+xq^{4j-3})}
{\prod_{j=1}^{\infty}(1-q^{4j})(1-x^2q^{4j-2})
\prod_{j=1}^{\infty}(1-q^{4j-2})}
\times\prod_{r=1}^{\infty}(1-q^{2r(m+1)})\\[5pt]
=&\frac{(-xq;q^2)_{\infty}(q^{2m+2};q^{2m+2})_{\infty}}{(q^2;q^2)_{\infty}
(x^2q^2;q^4)_{\infty}}.
\end{align*}
This completes the proof.  \qed

Similarly, one can check that
$$\Phi_{E,2m+2}(xq,xq,x^{-1}q,x^{-1}q)$$
is the generating function for partitions such that each even part appears
at most $2m+1$ times and the exponent of $x$ is
the alternating sum. Moreover, one sees that
$$\Psi_{E,2m+1}(xq,x^{-1}q,xq,x^{-1}q)$$
is the generating function for partitions such that
each even part appears at most $2m+1$ times and
the exponent of $x$ is the number of odd parts.
Thus Theorem \ref{partition theorem2}
is a consequence of the following relations.

\begin{thm}
We have
\begin{equation*}
\Phi_{E,2m+2}(xq,xq,x^{-1}q,x^{-1}q)
=\Psi_{E,2m+1}(xq,x^{-1}q,xq,x^{-1}q)
=\frac{(-xq;q^2)_{\infty}(q^{4m+4};q^{4m+4})_{\infty}}{(q^2;q^2)_{\infty}
(x^2q^2;q^4)_{\infty}}.\label{4eq2}
\end{equation*}
\end{thm}

{\noindent\it Proof.} Making the substitutions $a=xq,$
$b=xq$, $c=x^{-1}q$, $d=x^{-1}q$, $R=E$
(the set of even positive integers) and
$\rho=2m+2$ in Theorem \ref{3-1}, we find that
\begin{align*}
\lefteqn{\Phi_{E,2m+2}(xq,xq,x^{-1}q,x^{-1}q)}\\[5pt]
=&\prod_{j=1}^{\infty}
\frac{(1+(xq)^{\lceil\frac{j}{2}\rceil}
(xq)^{\lfloor\frac{j}{2}\rfloor}
(x^{-1}q)^{\lceil\frac{j-1}{2}\rceil}
(x^{-1}q)^{\lfloor\frac{j-1}{2}\rfloor})}
{(1-(xq)^{\lceil\frac{j}{2}\rceil}
(xq)^{\lfloor\frac{j}{2}\rfloor}
(x^{-1}q)^{\lceil\frac{j}{2}\rceil}
(x^{-1}q)^{\lfloor\frac{j}{2}\rfloor})
(1-(xq)^{j}(xq)^{j}(x^{-1}q)^{j-1}
(x^{-1}q)^{j-1})}\\[5pt]
&\times\prod_{r\in E}
(1-(xq)^{\lceil\frac{r}{2}\rceil
\frac{2m+2}{2}} (xq)^{\lfloor
\frac{r}{2}\rfloor\frac{2m+2}{2}}
(x^{-1}q)^{\lceil\frac{r}{2}\rceil
\frac{2m+2}{2}} (x^{-1}q)^{\lfloor
\frac{r}{2}\rfloor\frac{2m+2}{2}})\\[5pt]
=&\frac{\prod_{j=1}^{\infty}
(1+(xq)^{j}(xq)^{j}(x^{-1}q)^{j}(x^{-1}q)^{j-1})
\prod_{j=1}^{\infty}(1+(xq)^{j}
(xq)^{j-1}(x^{-1}q)^{j-1}(x^{-1}q)^{j-1})}
{\prod_{j=1}^{\infty}(1-(xq)^{\lceil\frac{j}{2}\rceil}
(xq)^{\lfloor\frac{j}{2}\rfloor}
(x^{-1}q)^{\lceil\frac{j}{2}\rceil}
(x^{-1}q)^{\lfloor\frac{j}{2}\rfloor})
\prod_{j=1}^{\infty}(1-x^2q^{4j-2})}\\[5pt]
&\times\prod_{r=1}^{\infty}(1-(xq)^{r(m+1)}
(xq)^{r(m+1)} (x^{-1}q)^{r(m+1)}(x^{-1}q)^{r(m+1)})\\[5pt]
=&\frac{\prod_{j=1}^{\infty}(1+xq^{4j-1})
\prod_{j=1}^{\infty}(1+xq^{4j-3})}
{\prod_{j=1}^{\infty}(1-q^{4j})(1-q^{4j-2})
\prod_{j=1}^{\infty}(1-x^2q^{4j-2})}
\times\prod_{r=1}^{\infty}(1-q^{4r(m+1)})\\[5pt]
=&\frac{(-xq;q^2)_{\infty}(q^{4m+4};q^{4m+4})_{\infty}}{(q^2;q^2)_{\infty}
(x^2q^2;q^4)_{\infty}}.
\end{align*}

On the other hand, making the  the
substitutions $a=xq$,
$b=x^{-1}q$, $R=E$  and
$\gamma=2m+1$ in  Theorem \ref{general two parameter},  we get
\begin{align*}
\lefteqn {\Psi_{E,2m+1}(xq,x^{-1}q,xq,x^{-1}q)}\\[5pt]
=&\prod_{j=1}^{\infty}
\frac{(1+(xq)^{\lceil\frac{j}{2}\rceil}
(x^{-1}q)^{\lfloor\frac{j}{2}\rfloor}
(xq)^{\lceil\frac{j-1}{2}\rceil}
(x^{-1}q)^{\lfloor\frac{j-1}{2}\rfloor})}
{(1-(xq)^{2\lceil\frac{j}{2}\rceil}
(x^{-1}q)^{2\lfloor\frac{j}{2}\rfloor})
(1-(xq)^{2j-1}(x^{-1}q)^{2j-1})}\\[5pt]
&\times\prod_{r\in E}
(1-(xq)^{\lceil\frac{r}{2}\rceil(2m+1+1)}
(x^{-1}q)^{\lfloor\frac{r}{2}\rfloor(2m+1+1)})\\[5pt]
=&\frac{\prod_{j=1}^{\infty}
(1+(xq)^{j}(x^{-1}q)^{j}(xq)^{j}(x^{-1}q)^{j-1})
\prod_{j=1}^{\infty}(1+(xq)^{j}(x^{-1}q)^{j-1}
(xq)^{j-1}(x^{-1}q)^{j-1})}
{\prod_{j=1}^{\infty}(1-(xq)^{2\lceil
\frac{j}{2}\rceil}(x^{-1}q)^{2\lfloor\frac{j}{2}\rfloor})
\prod_{j=1}^{\infty}(1-q^{4j-2})}\\[5pt]
&\times\prod_{r=1}^{\infty}(1-(xq)^{2r(m+1)}
(x^{-1}q)^{2r(m+1)})\\[5pt]
=&\frac{\prod_{j=1}^{\infty}(1+xq^{4j-1})
\prod_{j=1}^{\infty}(1+xq^{4j-3})}
{\prod_{j=1}^{\infty}(1-q^{4j})(1-x^2q^{4j-2})
\prod_{j=1}^{\infty}(1-q^{4j-2})}
\times\prod_{r=1}^{\infty}(1-q^{4r(m+1)})\\[5pt]
=&\frac{(-xq;q^2)_{\infty}(q^{4m+4};q^{4m+4})_{\infty}}{(q^2;q^2)_{\infty}
(x^2q^2;q^4)_{\infty}}.
\end{align*}
 This completes the proof. \qed

\section{Combinatorial Proof of Theorem 3.1}

In this section, we present a combinatorial
 proof of Theorem \ref{partition theorem1}.
 We need a property of
 Sylvester's bijection  between the set
$\mathcal{D}(n)=\{\lambda\mapsto n: \mbox{the parts of } \lambda \mbox{ are distinct}\}$
and the set $\mathcal{O}(n)=\{\lambda\mapsto n: \mbox{all parts of } \lambda \mbox{ are odd}\}$,
see also \cite{Pak06}.
Bessenrodt \cite{Bes94} showed that
Sylvester's bijection maps a partition $\lambda\in\mathcal{D}(n)$ to a partition
$\tau\in\mathcal{O}(n)$ such that
\begin{equation} \label{lambda1}
\lambda_1=l(\tau)+(\tau_1-1)/2
\end{equation}
and
the alternating sum of $\lambda$ equals the number of odd parts of $\tau$, that is,
\begin{equation}\label{lambda2}
l_a(\lambda)=l_o(\tau).
\end{equation}

We shall give a bijection $\Psi$ between the set $A_m(n)$
of partitions of $n$ such that each part appears at most $2m+1$ times
and the set $B_m(n)$ of partitions of $n$
such that each even part appears at most $m$ times.
Moreover, for any partition $\alpha$ in $A_m(n)$, we  show that
 the alternating sum of $\alpha$
is equal to the number of odd parts of $\Psi(\alpha)$ in $B_m(n)$.

Let $\alpha$ be a partition in $A_{m}(n)$. The map $\Psi$ from $A_m(n)$ to $B_m(n)$ can be
described as follows.

\begin{enumerate}
\item Write $\alpha$  as
a bipartition $\phi_1(\alpha)=(\lambda,\mu)$ subject to the following conditions.
First, set $\lambda$ and $\mu$ to be the empty partition.
For each part $t$ of $\alpha$, if $t$ appears an odd number of times, then
add one part $t$ to $\lambda$ and move the remaining parts $t$ to $\mu$; otherwise, move all parts $t$ of $\alpha$ to $\mu$.
Clearly,
$\lambda$ is a partition with distinct parts and each part in $\mu$ appears
an even number times. Moreover,
 each part of $\mu$ appears at most $2m$ times.

\item We apply  Sylvester's bijection $\phi_2$
to $\lambda$ to obtain
a partition $\phi_2(\lambda)=\tau$ containing only odd parts.

\item Since each part $t$ in $\mu$ appears
an even number times,  we merge two parts $t$ into a single part $2t$.
This leads to a partition $\phi_3(\mu)=\nu$ consisting of even parts such that
each even part appears at most $m$ times.

\item Putting the parts of $\tau$ and $\nu$ together, we obtain
a partition
$\phi_4(\tau,\nu)=\beta.$  It is clear that
$\beta\in B_{m}(n)$.
\end{enumerate}

Let $\beta$ be a partition in $B_{m}(n)$. The inverse map $\Psi^{-1}$
from $B_m(n)$ to $A_m(n)$ can be described as follows.
\begin{enumerate}
\item Write $\beta$ as a bipartition $\phi_4^{-1}(\beta)=(\tau,
\nu)$, where $\tau$ consists of the odd parts of $\beta$ and $\nu$ consists
of the even parts of $\beta$.
Since each even part of $\beta$
appears at most $m$ times,
the multiplicity of each part in $\nu$ does not exceed  $m$.

\item Decompose each even part $2t$ of $\nu$ into two equal
  parts $t$. We obtain a partition $\phi_3^{-1}(\nu)=\mu$
such that each part in $\mu$ appears an even number times.
Clearly,  each part of $\mu$ appears at most $2m$ times.

\item Since $\tau$ is a partition with odd parts, we can apply Sylvester's bijection
$\phi_2^{-1}$ to $\tau$ to obtain
a partition $\phi_2^{-1}(\tau)=\lambda$
with distinct parts.

\item Putting the parts of $\lambda$ and $\mu$ together,
we obtain a partition
$\phi_1^{-1}(\lambda,\mu)=\alpha$.
It is easy to see that $\alpha\in A_{m}(n)$.
\end{enumerate}

It is readily seen that the above map $\Psi$ is a bijection.
Here is an example for $m=2$ and $n=48$.
Let
$$\alpha=(1,2,2,2,2,2,4,4,4,4,7,7,7)\in A_2(48).$$
We decompose $\alpha$ into a bipartition
$$\phi_1(\alpha)=(\lambda,\, \mu)=
((1,2,7),(2,2,2,2,4,4,4,4,7,7)).$$
Then we apply Sylvester's bijection to
$\lambda$ to obtain
$$\phi_2(\lambda)=\tau=(1,1,1,1,3,3).$$
Merging equal parts in $\mu$, we get
$$\phi_3(\mu)=\nu=(4,4,8,8,14).$$
Finally, putting the parts of $\tau$ and $\mu$ together
we obtain $$\phi_4(\tau,\nu)=
\beta=(1,1,1,1,3,3,4,4,8,8,14)\in B_2(48).$$

The above bijection $\Psi$ leads to a combinatorial proof of
Theorem \ref{partition theorem1}.

{\noindent \it Proof of Theorem \ref{partition theorem1}.}
For a  partition $\alpha\in A_m(n)$, let $\beta=\Psi(\alpha) \in B_m(n)$.
We aim to show that the alternating sum of $\alpha$ is equal to
the number of odd parts in $\beta$.

For $\phi_1(\alpha)=(\lambda,\mu)$ and any part in $\mu$ appears
 an even number times in $\mu$, the alternating sum of $\alpha$ is equal to
the alternating sum of $\lambda$.
Using Bessenrodt's refinement \cite{Bes94},
it is easy to see that the
alternating sum of $\lambda$ is equal to the
number of odd parts in
$\tau$. Since all the parts in $\nu$ are even numbers, the
number of odd parts in $\tau$ is just the number of odd
parts of $\beta.$ Consequently, we see that the alternating sum of $\alpha$ is equal to
the number of odd parts in $\beta$. Hence Theorem \ref{partition theorem1} holds. \qed

We conclude with  a stronger version of
 Theorem \ref{partition theorem1} based on
the above bijection $\Psi$.
Let $\mathcal{A}_{\varphi}(n)$ be the set of partitions of $n$ such that the number of
appearance of part $i$ is at most $2\varphi(i)+1$, for any $i>0$.
Let $\mathcal{B}_{\varphi}(n)$ be the set of partitions of $n$ such that the number of
appearance of part $2i$ is at most $\varphi(i)$, for any $i>0$.
Notice that for any map $\varphi$ from  $\mathbb{N}$
to $\mathbb{N}$, the bijection $\Psi$ can be applied to  $\mathcal{A}_{\varphi}(n)$
and  $\mathcal{B}_{\varphi}(n)$.

Let ${p}(\lambda)$ denote the largest odd part
in $\lambda$, and let $q(\lambda)$   denote the largest
part of $\lambda$ with odd multiplicity. Then we have the following correspondence.

\begin{thm}\label{general of partition theorem1}
Let $P(n,\varphi,i,k,t)$ be the set of partitions
 $\alpha\in\mathcal{A}_{\varphi}(n)$ of $n$ such that the number of
appearance of $i$ in $\alpha$ is at
most $2\varphi(i)+1$, $l_a(\alpha)=k$, and $q(\alpha)=t$.
Let $Q(n,\varphi,i,k,t)$ be the set of partitions
$\beta\in\mathcal{B}_{\varphi}(n)$ of $n$ such that
the number of appearance of $2i$ in $\beta$ is
 at most $\varphi(i)$, $l_o(\beta)=k$,
 and $l_o(\beta)+(p(\beta)-1)/2=t$. Then for $i\ge 1$ and $k,t\ge 0$, there
 exists a bijection between $P(n,\varphi,i,k,t)$ and $Q(n,\varphi,i,k,t)$.
\end{thm}

{\noindent \it Proof of Theorem
\ref{general of partition theorem1}.} Let $\alpha\in P(n,\varphi,i,k,t)$,
that is, the number of appearances
of $i$ in $\alpha$ is at most $2\varphi(i)+1$, $l_a(\alpha)=k$,
and $q(\alpha)=t$. Let
$\Psi$ be the map from $\mathcal{A}_{\varphi}(n)$ to
$\mathcal{B}_{\varphi}(n)$.
Denote $\Psi(\alpha)$ by $\beta$. We aim to show that
$\beta\in Q(n,\varphi,i,k,t)$, that is, the number of appearances
of $2i$ in $\beta$ is at most $\varphi(i)$, $l_o(\beta)=k$,
and $l_o(\beta)+(p(\beta)-1)/2=t$.

First, write $\alpha$  as
a bipartition $\phi_1(\alpha)=(\lambda,\mu)$ based on the following procedure.
Set the initial values of $\lambda$ and $\mu$ to be the empty partition.
For each part $t$ of $\alpha$, if $t$ appears an odd number of times in $\alpha$, then
add a part $t$ to $\lambda$ and move the remaining parts $t$ of $\alpha$
to $\mu$; otherwise, move all parts $t$ of $\alpha$ to $\mu$. Now, $\lambda$ is a partition with distinct parts with
$\lambda_1=q(\alpha)=t$. Moreover, since any part $i$ of $\alpha$ appears at most
$2\varphi(i)+1$ times,  the part $i$ in $\mu$ appears at most
$2\varphi(i)$ times. It is evident that $l_a(\lambda)=\l_a(\alpha)=k$ .

Since $\lambda$ is a partition with distinct parts, we can apply Sylvester's
bijection $\phi_2$ to $\lambda$ to obtain a partition of odd parts. Let
$\tau=\phi_2(\lambda)$. Recall that  $l(\tau)+(\tau_1-1)/2=\lambda_1=t$ and
$l(\tau)=\l_a(\lambda)=k$, as given in  (\ref{lambda1}) and (\ref{lambda2}).
Furthermore, we can merge every two equal parts in $\mu$ to form
a partition with even parts, denoted $\nu=\phi_3(\mu)$. Observe that
 the number of
 appearance of part $2i$ in $\nu$ is at most
 $\varphi(i)$.

 Finally, we put the parts of $\tau$ and $\nu$ together to form
a partition $\beta=\phi_4(\tau,\nu)$.
It can be seen that the number of appearance of part $2i$ in
$\beta$ is at most
 $\varphi(i)$ and $l_o(\beta)=l(\tau)=k$. Since $\nu$ is partition with even parts,
it is easy to see that $p(\beta)=\tau_1$.
Using the facts $l_o(\beta)=l(\tau)$, $p(\beta)=\tau_1$
and $l(\tau)+(\tau_1-1)/2=t$, we obtain that
$l_o(\beta)+(p(\beta)-1)/2=t$.
This completes the proof. \qed

\section{Combinatorial Proof of Theorem 3.2}

In this section, we prove Theorem~3.2 combinatorially.   We perform the same separation procedure $\phi_1$ described in Step 1 of the bijection $\Psi$ on a partition $\alpha \in B_{2m+1}(n)$ to obtain a pair $(\lambda, \mu)$, where $\lambda$ is a partition into distinct parts and $\mu$ is a partition into parts with even multiplicities. We then apply the Sylvester bijection $\phi_2$ described in Step 2 of $\Psi$ to $\lambda$ and denote the resulting partition by $\tau$. For $\mu$, we take each of the odd parts $2i-1$ and write its multiplicity $m_{2i-1}$, which is even,  as the sum of powers of $2$, namely
\begin{align*}
m_{2i-1}= a_1 2^1+a_2 2^2+\cdots,
\end{align*}
where $a_j=0$ or $1$. We map the copies of odd parts $2i-1$ from $\mu$ to $2^j(2i-1)$ with multiplicity $a_j$ for $j\ge 1$. Since $a_j=0$ or $1$, at most one more copy of each even part can be obtained. Combining these with the even parts of $\mu$, we obtain a partition $\nu$ into even parts with multiplicity at most $2m+1$.  We define $\beta$ to be the partition whose parts are the union of the parts of $\tau$ and $\nu$.

We now show that the map defined above is indeed a bijection. Since the only difference from the proof of Theorem~3.1 is the map for $\mu$, it is sufficient to show that this map is reversible.  Note that every even integer can be uniquely written as $2^j(2i-1)$ for some positive integers  $i$ and $j$. Thus, given a partition $\nu$ into even parts $2^j(2i-1)$ with multiplicities $m_{2^j(2i-1)}$ not exceeding $2m+1$, if $m_{2^j(2i-1)}$ is odd, then we transform one copy of the part $2^j(2i-1)$ into $2^j$ copies of the part $2i-1$. If $m_{2^j(2i-1)}$ is even, then it is less than $2m+1$ and we do nothing. After transforming one copy of each part of $\nu$ with an odd multiplicity, we obtain a partition $\mu$, where each  part  occurs an even number of times, in particular its even parts cannot occur more than $2m$ times.  This completes the proof.

 \vspace{.3cm}

\noindent{\bf Acknowledgments.} This work was
supported by the 973 Project,
the PCSIRT Project of the Ministry of Education,
and the National Science Foundation of China.


\vskip 10pt

\begin{thebibliography}{99} \small
\setlength{\itemsep}{-.8mm}

\bibitem{Ald69} H.L. Alder, Partition
identities -- From Euler to the present,
 Amer. Math. Monthly 76 (1969), 733--746.

\bibitem{And66} G.E. Andrews, On generalizations
 of Euler's partition theorem, Michigan Math. J. 1 (1966),
 491--498.

\bibitem{And72} G.E. Andrews, Partition identities,
 Adv. Math. 9 (1972), 10--51.

\bibitem{And76} G.E. Andrews,
The Theory of Partitions,
Addison-Wesley Publishing Co., 1976.

\bibitem{AndE04} G.E. Andrews and Kimmo Eriksson,
Integer Partitions, Cambridge Unversity Press, 2004.


\bibitem{Bes94} C. Bessenrodt, A bijection for
Lebesgue's partition identity in the
spirit of Sylvester, Discrete Math. 132 (1994), 1--10.

\bibitem{Bou971} M. Bousquet-M\'{e}lou and K. Eriksson,
Lecture hall partitions, Ramanujan J.  1 (1997), 101--111.

\bibitem{Bou972}M. Bousquet-M\'{e}lou and K. Eriksson,
Lecture hall partitions II, Ramanujan J.  1 (1997), 165--185.

\bibitem{Bou06} C. Boulet,
 A four parameter partition identity,
 Ramanujan J. 12(3) (2006), 315--320.

\bibitem{Chen10} W.Y.C. Chen, H.Y. Gao,
K.Q. Ji and M.Y.X. Li, A unification of
two refinements of Euler's partition
theorem, Ramanujan J. 23 (1) (2010), 137-149.

\bibitem{Fin88} N.J. Fine, Basic Hypergeometric
 Series and Applications, Math. Surveys
27, AMS Providence, 1988.

\bibitem{Gla83} J.W.L. Glaisher,
A theorem in partitions, Messenger of Math.
12 (1883), 158--170.

\bibitem{Gas90} G. Gasper and M. Rahman,
 Basic Hypergeometric Series, Encycl. of
Math. and Its Appl., Vol. 35,
 Cambridge University Press, Cambridge,
1990.

\bibitem{Moo74} E. Moore,
Generalized Euler-type partition
identities, J. Combin. Theory Ser.
A 17 (1974), 78--83.

\bibitem{Pak06} I. Pak,
Partition bijections, a survey,
Ramanujan J. 12(1) (2006), 5--75.

\bibitem{Syl82} J. Sylvester,
A constructive theory of partitions,
arranged in three acts, an interact and
an exodion, Amer. J. Math. 5 (1882), 251--330.

\end{thebibliography}
\end{document}